\documentclass[11pt,a4paper]{amsart}

\usepackage{amsmath,amsfonts,amsthm,amssymb}

\title{On Morrison's cone conjecture\\ 
for klt surfaces with $K_X\equiv0$}
\author{Kaori Suzuki}
\date{}
\address{Department of Mathematical Sciences, University of
Tokyo, Komaba, Meguro, Tokyo 153--8914, Japan}
\email{suzuki@ms.u-tokyo.ac.jp}

\theoremstyle{definition}
\newtheorem{definition}{Definition}[section]

\theoremstyle{plain}
\newtheorem{thm}[definition]{Theorem}
\newtheorem{prop}[definition]{Proposition}
\newtheorem{lemma}[definition]{Lemma}

\newtheorem{conj}[definition]{Conjecture}
\theoremstyle{remark}
\newtheorem{rem}[definition]{Remark}
\newtheorem{cl}[definition]{Claim}

\newtheoremstyle{citing}
{3pt}{3pt}{\itshape}{}{\bfseries}{.}{.5em}{\thmnote{#3}}

\theoremstyle{citing}

 \newcommand{\sA}{\mathcal A}
 \newcommand{\sC}{\mathcal C}
 \newcommand{\Oh}{\mathcal O}
 \newcommand{\C}{\mathbb C}
 \newcommand{\Q}{\mathbb Q}
 \newcommand{\R}{\mathbb R}
 \newcommand{\Z}{\mathbb Z}
 \newcommand{\directsum}{\bigoplus}
 \newcommand{\iso}{\simeq}
 \newcommand{\tensor}{\otimes}
 \newcommand{\id}{\mathrm{id}}
 \newcommand{\ga}{\gamma}
 \newcommand{\Ga}{\Gamma}
 \newcommand{\De}{\Delta}
 \newcommand{\fie}{\varphi}
 \newcommand{\si}{\sigma}
 \newcommand{\om}{\omega}
\DeclareMathOperator{\Ann}{Ann}
\DeclareMathOperator{\Aut}{Aut}
\DeclareMathOperator{\diag}{diag}
\DeclareMathOperator{\GL}{GL}
\DeclareMathOperator{\Hom}{Hom}
\DeclareMathOperator{\im}{im}
\DeclareMathOperator{\Int}{Int}
\DeclareMathOperator{\ord}{ord}
\DeclareMathOperator{\pr}{pr}

\DeclareMathOperator{\rank}{rank}
\DeclareMathOperator{\Spec}{Spec}

\newcommand{\Span}[1]{\left<#1\right>}
\newcommand{\rest}[1]{{}_{{\textstyle|}#1}}
\DeclareMathOperator{\ns}{\mathit{N}^1}

 \def\coh2#1{{\mathit{H}^2({#1},\Z)}{}}
\numberwithin{equation}{section}
\renewcommand{\labelenumi}{\theenumi.}

\begin{document}

\begin{abstract} This paper considers normal projective complex surfaces
$X$ with at worst Kawamata log terminal singularities and $K_X\equiv0$.
The aim is to prove that there is a finite rational polyhedral cone which
is a fundamental domain for the action $\Aut X$ on the convex hull of its
rational ample cone.
\end{abstract}

\maketitle

\renewcommand{\rightmark}{Morrison's cone conjecture for surfaces with
$K\equiv0$}

\setcounter{section}{-1}
\section{Introduction}\label{sec!In}

Let $X$ be a normal projective complex surface with at worst Kawamata log
terminal singularities (klt singularities) and $K_X \equiv0$. By the
classification of surfaces, $X$ is either an Abelian surface, a
hyperelliptic surface, a K3 surface with only rational double points
(RDPs), an Enriques surface with only RDPs or a log Enriques surface. Here
a {\em log Enriques surface} (compare Zhang \cite{Zhang}) is a rational
normal projective surface with at worst klt singularities and $K_X\equiv0$.

The aim of this paper is to study the $2$-dimensional klt analogue of the
following conjecture of Morrison:

\begin{conj}[Morrison's Cone Conjecture {\cite[Section 4]{morrison3}}]
Let $V$ be a Calabi--Yau manifold, $\sA(V)$ its ample cone, and $\sA'(V)$
the convex hull of its rational ample cone $\overline{\sA(V)}\cap
H^2(V,\Q)$. Then there exists a rational polyhedral cone $\De
\subset \sA'(V)$ such that $\Aut V\cdot\De=\sA'(V)$.
\end{conj}

Our main result is as follows:

\begin{thm}[Main Theorem]\label{main theorem}\label{th!main} Let $X$ be a
normal projective complex surface with at worst klt singularities and
$K_X\equiv0$. Then there exists a rational finite polyhedral cone $\De$
which is a fundamental domain for the action of $\Aut X$ on $\sA'(X)$ in
the sense that
\begin{enumerate}
 \item $\sA'(X)=\bigcup_{\theta\in \Aut X}\theta_*\De$,
 \item $\Int\De \cap\theta_*\Int\De=\emptyset$ 
unless $\theta_*=\id$.
\end{enumerate}
\end{thm}

This conjecture is known for a smooth Abelian variety (Kawamata \cite{k2})
or a smooth K3 surface (Sterk \cite{sterk}). Our new result covers the
remaining cases.

The plan of the proof is as follows: for $X$ as in Theorem~\ref{th!main},
write $I=I(X)$ for the (global) index of $X$, that is, the smallest
positive integer such that $IK_X$ is Cartier and linearly equivalent to
zero and $\pi\colon Y:=\Spec(\directsum_{i=0}^{I-1}\Oh_X(-iK_X))\to X$ for
the (global) index one cover of $X$. Then $\pi\colon Y\to X$ is a Galois
finite cover with Galois group $G=\Span g\iso\Z/I$, and
$g^*(\om_Y)=\zeta_I\om_Y$, where $\om_Y$ is a generator of
$H^0(\Oh_Y(K_Y))$ and $\zeta_I$ a primitive $I$th root of 1. It is known
that $Y\to X$ is etale in codimension~1, so that $G$ admits no pointwise
fixed curve. Since $Y$ is a projective surface with only RDPs and
$\Oh_Y(K_Y)\iso\Oh_Y$, it is either a projective K3 surface with only RDPs
or a smooth Abelian surface.  We treat log Enriques surfaces according to
this case division.

In Section~1, we consider the case where the global index one cover $Y$ is
an Abelian surface. In this case the theorem follows from Kawamata's result
\cite{k2} and the Torelli theorem for Abelian surface (Shioda
\cite{shioda}). In the course of the proof, we also classify log Enriques
surfaces whose index one cover is an Abelian surface. Section~2 considers
the case where $Y$ is a K3 surface with at worst RDPs. In this case, we use
the equivariant version of Torelli theorem (Oguiso and Sakurai
\cite{sakurai}) and Morrison's weakly polarized global Torelli theorem for
K3 surfaces \cite{M}. Our result Theorem~\ref{th!K3} is slightly stronger
than we need.

\subsection*{Acknowledgements} The author would like to express her deep
gratitude to Professor Keiji Oguiso for his valuable comments and
unceasing encouragement. She also thanks Professor Yujiro Kawamata who for
stimulating discussion and Professor Miles Reid for help with the English.

\section{The case where the index one cover is an Abelian surface}

Assume first that $X$ is an Abelian surface. In this case, we obtain the
result at once from the following theorem. 

 \begin{thm}[Kawamata {\cite[Theorem~2.1]{k2}}]\label{kawamata} Assume
that $X$ is an\linebreak Abelian surface. Then, there exists a finite
rational polyhedral cone $\De$ which is a fundamental domain for the
action of\/ $\Aut X$ on $\sA'(X)$ in the sense that
\begin{enumerate}
 \item $\sA'(X)=\bigcup_{\theta\in \Aut X}\theta_*\De$,
 \item $\Int\De \cap\theta_*\Int\De=\emptyset$ 
unless $\theta_*=\id$.
\end{enumerate}
\end{thm}

The case of a hyperelliptic surface $X$ is trivial: Since $\rho(X)=2$ and $\bar{NS}X$ is spanned by its two 
elliptic fibrations, each of which is invariant under $\Aut X$.

Next, consider the case of a log Enriques surface $X$ whose index one
cover $Y$ is a smooth Abelian surface.

\begin{prop}\label{prop1}
Let $X$ be a log Enriques surface whose index one cover $(Y, \Span g)$ is
an Abelian surface. Then $Y \to X$ is isomorphic to one of the following
two cases:
\begin{enumerate}
\item $Y=E_{\zeta_3}^2$ and $g_3=\diag(\zeta_3,\zeta_3)$; here
$E_{\zeta_3}$ is the elliptic curve with multiplication by $\zeta_3$. Or
\item $Y=$ the Jacobian surface $J(C)$ of the genus two curve
$C:y^2=x^5-1$, and the automorphism $g_5$ of order $5$ is induced by the
automorphism $(x,y)\mapsto(\zeta_5x,y)$ of $C$.
\end{enumerate}
\end{prop}

Our proof of Proposition~\ref{prop1} is similar to Oguiso
\cite[Lemma~1.1]{type3}.

\begin{lemma}[Oguiso]\label{lemma1.1} Let $X$ be a log Enriques surface
whose index one cover $(Y,\Span g)$ is an Abelian surface. Then
\begin{enumerate}
 \item $I=\ord g=3$ or $5$;
 \item $(Y,\Span g)$ is isomorphic to $(E_{\zeta_3}^2,\Span{g_3})$ if $I=3$;
 \item $(Y,\Span g)$ is isomorphic to $(J(C),\Span{g_5})$ if
$I=5$.
\end{enumerate} 
\end{lemma}

\begin{proof} Write $Y$ as $V/\Lambda$, where $V=\C^2$ and $\Lambda$ is a
discrete $\Z$-submodule of $V$ of rank $4$. Since $G$ does not have any
nonzero translations, $\Span g$ has a fixed point:
\begin{cl}
$g$ has a fixed point, that is, there exists $p\in Y$ 
such that $g(p)=p$.
\end{cl}
\begin{proof} Since $\Aut Y=Y\rtimes \Aut_{\mathrm{Lie}}(Y)$, in suitable
coordinates, we can write
 \[
g^*\begin{pmatrix}
x \\ y
\end{pmatrix}=A
\begin{pmatrix}
x \\ y
\end{pmatrix}+\begin{pmatrix}
a_1 \\ a_2
\end{pmatrix}.
 \]
Since $g$ is of finite order, we may take $A$ as a diagonal matrix 
$\begin{pmatrix} \alpha & 0 \\ 0 & \beta \end{pmatrix}$.

It suffices to show that $\alpha,\beta\ne1$. We derive a contradition from
the assumption $\alpha=1$. Recall that $Y^{\Span g}\ne\emptyset$. Then
there exists $k<\ord g$ such that $Y^{g^k}\ne\emptyset$. We can write
 \[
{g^k}^*
\begin{pmatrix}
x \\ y
\end{pmatrix}=
\begin{pmatrix}
1 & 0 \\
0 & \beta^k
\end{pmatrix}
\begin{pmatrix}
x \\ y
\end{pmatrix}+\begin{pmatrix}
b_1 \\ b_2
\end{pmatrix} \quad \hbox{for some $b_1, b_2$.}
 \]
Let $p=\begin{pmatrix} x \\ y \end{pmatrix}\in Y^{g^k}$. Then $x=x+b_1$
leads to $b_1=0$. Thus the curve $(y=0)$ in $Y$ is contained in $Y^{\Span
g}$. This contradicts the fact that $Y^{\Span g}$ is a finite set.
\end{proof}

Since $Y$ is an Abelian surface, there exist global coordinates $(x,y)$ of
$Y$ centered at a fixed point of $g$ such that $g^*(x,y)=(ax,by)$. Recall
that $g^*\om_Y=\zeta_I\om_Y$, where $\ord g=I$ and $Y^{\Span g}$ is
a finite set. Then we have $g=\diag(\zeta_I^{k_1},\zeta_I^{k_2})$ for some
integers $k_1$ and $k_2$ where each $k_i$ is coprime to $I$. 
 
Let $\Z[\Span g]$ be the group algebra of the cyclic group $\Span g$ over
$\Z$. Clearly, $\Z[\Span g]\iso\Z[x]/(x^I-1)$. The action of $g$ on $V$
defines a natural $\Z[\Span g]$-module structure on $\Lambda$:
$F(g)(\lambda)=\diag(F(\zeta_I^{k_1}),F(\zeta_I^{k_2}))\cdot \lambda$.

Write $\Phi_I(x)$ for the $I$th cyclotomic polynomial. Since
$\Ann\Lambda=(\Phi_I(g))$ and $\Lambda$ is a torsion free $\Z[\Span
g]/(\Phi_I(\Span g)$-module. This also holds for $V$.

Tensoring with $\Q$ the tower of modules
$\Z\subset\Z[\zeta_I]\subset\Lambda$, gives the tower of vector spaces
$\Q\subset\Q[\zeta_I]\subset\Lambda\tensor\Q$. Then we get
$\phi(d)\cdot\dim_{\Q[\zeta_I]}(\Lambda\tensor_{\Z}\Q)=4$, where $\phi$ is
the Euler function. Thus, $\phi(d)$ is either 2 or 4. Combining this with
the fact $I$ is not divisible by 2 (\cite[Lemma~2.2]{Zhang}), we get $I$
is either 3 or 5.

% Let $\lambda$ be an element of either $\Lambda$ or $V$ and 

%%%%%%%%%%%%%%%%%%%%%%%%%%%%% I=3 case 
For the case of $I=3$, we have $g^*\om_Y=\zeta_3\om_Y$. Therefore
$ab=\zeta_3$. We also see that $a$ and $b$ must be roots of unity. Recall
that neither $a$ nor $b$ is $1$. Thus, we conclude that
$(a,b)=(\zeta_3^2,\zeta_3^2)$. Changing the generator by $g$ to $g^2$, we
obtain $g^*(x, y)=(\zeta_3x,\zeta_3y)$. Since $V$ is the direct sum of
some eigenspaces corresponding to the indicated eigenvalues of the action
of $\zeta_I \tensor 1$ on the space $\Lambda \tensor_{\Z}\Z$, $\Lambda$
must be a free $\Z[\zeta_I]$-module of rank $4/\phi(I)$. This gives
$\phi_3\colon\Lambda\iso\Z[\zeta_3]^{\oplus2}$.

 Next we consider the eigenspace decomposition:
 \[
\Z[\zeta_I] \tensor_{\Z} \C= \bigoplus_{\begin{smallmatrix} 0 < i < m\\
(m,d)=1 \end{smallmatrix}} V(m),
 \]
where $V(m)$ denotes the eigenspace of the eigenvalue $\zeta_I^m$. In our
case, we obtain the following $\phi_I$-equivariant eigenspace
decompositions: 
 \[
 \Z[\zeta_3]^{\oplus2} \tensor_{\Z} \C=V(1)^{\oplus2}\oplus
V(2)^{\oplus2}
 \quad\hbox{and}\quad
\Lambda \tensor_{\Z}\C=V'(1)\oplus V'(2).
 \]
 
Therefore we have $\phi_3\colon (Y,g)=(V'(1)/\Lambda,\zeta_3)\iso
(V(1)^{\oplus 2}/\Z[\zeta_3],\zeta_3)$. In particular, $(E_{\zeta_3}^2,
g_3)\iso (V(1)^{\oplus 2}/\Z[\zeta_3],\zeta_3)$.
 
%%%%%%%%%%%%%%%%%%%%%%%%%%%%% I=5 case
If $I=5$, then the same argument gives $(a, b)=(\zeta_5^3,\zeta_5^3)$, 
$(\zeta_5^2,\zeta_5^4)$ or $(\zeta_5^4,\zeta_5^2)$. Changing the
generator again, we obtain $g^*(x, y)=(\zeta_5x,\zeta_5y)$,
$(\zeta_5x,\zeta_5^2y)$. However the first case cannot happen because
$\rank\Lambda=4$. Therefore, $g^*(x, y)=(\zeta_5x,\zeta_5^2y)$. Then in a
similar way we have $\phi_5\colon\Lambda\iso\Z[\zeta_5]$.

Since $\Z[\zeta_5] \tensor_{\Z} \C=\directsum_{m=1}^4 V(m)$ and $\Lambda
\tensor_{\Z}\C=\directsum_{m=1}^4 V'(m)$, we have $\phi_5\colon
(Y,g)=(V(1)\oplus V(2))/\Z[\zeta_5],\zeta_5)$ and $(J(C),g_5)\iso
(V(1)\oplus V(2))/\Z[\zeta_5],\zeta_5)$. This gives the result.
\end{proof}

%%%%%%%%%%%%%%%%%%%%% log Enriques case
In the rest of the proof $X$ is a log Enriques surface; the proof reduces
to solving the following proposition, with $Y \to X$ the index one cover of
$X$ and $g$ a generator of the Galois group:

\begin{prop}\label{hifuku}
There exists a finite rational polyhedral cone $\De'$ inside $\sA'(Y)$
such that $\sA'(Y)^{\Span g}=\Aut(Y,g)\De'$.
\end{prop}

We need the following theorem due to Shioda:

 \begin{thm}[{\cite[Theorem~1]{shioda}}]\label{shioda1} Let $Y$ and $Y'$
be Abelian surfaces over $\C$, and let
 \[
 \fie\colon \coh2{Y}\to \coh2{Y'}
 \]
be an isomorphism of Euclidean lattices. 
Denote the period map of $Y$ (resp. Y') by 
$\mathit{p}_Y$ (resp. $\mathit{p}_{Y'}$).  
Then the necessary and sufficient
condition for the isomorphism $\fie$ or $-\fie$ to be induced by an
isomorphism $f\colon Y'\to Y$ is that $\fie$ satisfies the
following conditions:
\begin{enumerate}
 \item $\det\fie=+1$, and 
 \item $p_Y'\circ\fie=\mathrm{const.} p_Y$ (that is, $\fie$
preserves the periods).
\end{enumerate}
\end{thm}

If $I=3$, we take a basis of $\ns{Y}$ as $D_i=\pr_i^*(0)$ for $i=1, 2$,
$E_1= \De$ and $E_2= \Ga_{f}$, where $\pr_i$ are the projections, $\De$
denotes the diagonal and $\Ga_{f}$ is the graph of $f\colon E_{\zeta_3}\to
E_{\zeta_3}$ induced by $p \mapsto\zeta_3 p$. Then the action of $g$ on
$\ns{Y}$ is trivial and Theorem~\ref{kawamata} and Theorem~\ref{shioda1}
give the result immediately.

In the case $I=5$, from $g^*\om=\zeta_5\om$, the eigenvalues of $g$ on
$T_Y$ are $\zeta_5,\zeta_5^2,\zeta_5^3$ and $\zeta_5^4$, and $T$ is the
minimal lattice which contains 2-form $\om$. Then we have $\rank T_Y=4$.
Since $\ns{Y}=T_Y^{\perp}$, we have $\rank \ns{Y}=2$ and the action of
$g$ on $\ns{Y}$ is trivial. Thus Theorem~\ref{kawamata} and
Theorem~\ref{shioda1} again give the result.

%%%%%%%%%%%%%%%%%%%%% "핢,B*ƒ(BA[ƒxƒ‹
\begin{rem}
By the reduction theory and explicit calculation, in the case $I=3$ we have
\begin{multline*}
 \renewcommand{\arraystretch}{1.5}
\R_{\ge 0} 
 \begin{pmatrix} 
 0 & 0 \\
 0 & 1
 \end{pmatrix}
 +\R_{\ge 0} 
 \begin{pmatrix} 
 1 & 0 \\
 0 & 1
 \end{pmatrix} \\
 \renewcommand{\arraystretch}{1.5}
 +
 \R_{\ge 0} 
 \begin{pmatrix} 
 2 & 1+\zeta_3 \\
 \overline{1+\zeta_3} & 2
 \end{pmatrix}
 +\R_{\ge 0} 
 \begin{pmatrix} 
 2 & 1-\zeta_3 \\
 \overline{1-\zeta_3} & 2
 \end{pmatrix}
\end{multline*}
 as $\De$, under the identification
 \[
 \im(\si\colon \Aut Y\to \GL(\ns{Y},\Z))\iso \GL(2,\Z[\zeta_3]).
 \].
\end{rem}
%%%%%%%%%%%%%%%%%%%% "핢,B*'(B‹'R‹È-Ê

\section{The case where the index one cover is a K3 surface \\
with only RDPs}

In this section, we prove the following Theorem~\ref{th!K3}, which is
slightly stronger than necessary for the main result. Indeed, we obtain
Theorem~\ref{th!main} as a corollary on taking $G$ as the Galois group of
the index one cover and applying Theorem~\ref{th!K3}.

Let $\sA'(Y)$ be the convex hull of the rational ample cone of $Y$, and
$\sA(Y)^G$ be its $G$-fixed part. $\sA(Y,G):=\bigl\{\phi\in Y\bigm|\phi
\circ g= g\circ \phi\hbox{ for any }g\in G\bigr\}$.

\begin{thm}\label{th!K3}
Let $Y$ be a K3 surface with only RDPs and $G$ a finite subgroup of $\Aut
Y$. Then there exists a finite rational polyhedral cone $\De'$ which
is a fundamental domain for the action of $\Aut(Y,G)$ on $\sA'(Y)^G$ in
the sense that
\begin{enumerate}
 \item $\sA'(Y)^G=\bigcup_{\theta\in \Aut(Y,G)}\theta_*\De'$,
 \item $\Int\De'\cap\theta_*\Int\De'=\emptyset$ 
unless $\theta_*=\id$.
\end{enumerate}
\end{thm}

Let $\nu\colon Z\to Y$ be the minimal resolution of $Y$ and
$\bigcup_{i=1}^mE_i$ the exceptional divisor. From now, we mainly work on
$Z$. 

Our proof proceeds on similar lines to Oguiso and Sakurai \cite{sakurai}.
We first fix the notation adopted in our argument below.

From the uniqueness of the minimal resolution, we have the natural
auto\-morphism of $Z$ which corresponds to $g$. We denote this by $g$
again. Let $\coh2{Z}^W=\bigl\{D\in \coh2{Z} \bigm|(D\cdot E_i)=0 \hbox{
for all $i$}\bigr\}$. We set
 \begin{align*} S'&:=H^{1,1}(Z,\Z)\cap \coh2{Z}^W; \\ S&:=(S')^G=
\bigl\{D\in S' \bigm| g^*(D)=D \hbox{ for all } g\in G
\bigr\};
\\ T&:=\bigl\{D'\in \coh2{Z}^W\bigm|(D'\cdot D)=0 \hbox{ for all } D\in S
\bigr\};
\\ S^*&:=\Hom (S,\Z); \\
\sA&:= (\sA')^G=\sA' \cap S_{\R}; \\
 \end{align*}
$\sC^{\circ}$ is the positive cone of $X$, that is, the connected component
of the space $\bigl\{D\in \ns{X}_{\R}\bigm| (D\cdot D)>0\bigr\}$
containing the ample classes.

$\sC'$ is the union of $\sC^{\circ}$ and all $\Q$-rational rays in the
boundary $\partial(\sC')^{\circ}$ of $\sC^{\circ}$. 

$\sC:=(\sC')^G$;

$Q_Z:=\bigl\{f\in \Aut Z \bigm| f\circ g=g\circ f \hbox{ for all }g\in
G\bigr\}$;

$O(S)^{\circ}:=$ the orthogonal group of the lattice $S'$ which preserves
the positive cone $C$.

$O(S)^+:=$ the set of $\si\in O(S)^{\circ}$ that are Hodge isometries of 
$\coh2{Z}^W$ and extend to Hodge isometries of $\coh2{Z}$ satisfying
$\si\circ g^*=g^* \circ\si$ for all $g\in G$ (the last condition implies
$\si(S)=S$);

$P(S):=\bigl\{\si\in O(S)^{\circ} \bigm|\si(A)=A \}$;
 \[
 P(S)^+:=
 \left\{ \si\in O(S)^{\circ} \left|
 \begin{array}{l} \si=f^* \rest{S} \hbox{ for some }f\in Q_Z \\
 \hbox{and $f$ comes from $f'\in Q_Y$} \end{array} \right\}\right..
 \]
In other words, $P(S)^+$ is the image of the group homomorphism of $Q_Z\to
P(S)$ which maps from $f$ to $f^*\rest{S}$; this implies that
$\si(\sA)=\sA$ for each $\tau\in P(S)^+$.

Then we have the following:

 \begin{lemma}\label{os1.3} 
 \begin{enumerate} \renewcommand{\labelenumi}{(\arabic{enumi})}
 \item $S$ is an even hyperbolic lattice if\/ $\rank S \ne 1$. 
 \item The interior $\sA^{\circ}$ of $\sA$ consists of the
$G$-invariant ample classes of $Y$ and is nonempty. 
 Moreover, $\sA^{\circ}=(\sA')^{\circ} \cap S_{\R}$.
 \item Set $O(S)^{++}:= \bigl\{\tau\in O(S)^{\circ} \bigm| \tau
\rest{S^*/S}=\id\bigr\}$. Then $O(S)^{++} \subset O(S)^+ \subset
O(S)^{\circ}$ with each inclusion of finite index.
 \end{enumerate}
 \end{lemma}

%%%%%%%%%%%%%%%%%%%%%%%%%% Lemma~4.3 ,BL(BØ-3/4 

 \begin{proof} 
 The assertions (1), (2) are clear.  The last statement (3) comes from the
fact that
 \[
 O(S)^{++}=\ker\bigl[O(S)^{\circ}\to \Aut(S^*/S)\bigl]
 \]
 and $\Aut(S^*/S)$ is a finite group. We need to show that
$O(S)^{++}\subset O(S)^+$, that is, if $\tau\in O(S)^{++}$ then $\tau\in
O(S)^+$. Since $\tau\rest{(S^*/S)}=\id$, $(\tau,\id)\in O(S)\times O(T)$.
This is a Hodge isometry of $\coh2{Z}$. Set $\tau\rest{\coh2{Z}^W}=\si$.
Then
$\si$ preserves the Hodge structure and satisfies
$(\si(z)\cdot\si(w))=(z\cdot w)$ for $z,w\in\coh2{Z}^W$. Therefore $\si$
is a Hodge isometry of $\coh2{Z}^W$. Moreover, since $\tau \circ
g^*=g^*\circ \tau$ on $S$, and $\si\circ g^*=g^*\circ\si=g^*$ on $T$, we
have $\si\circ g^*=g^* \circ\si=\tau$ on $\coh2{Z}^W$. This shows $\tau\in
O(S)^+$.
 \end{proof}

Next we consider the reflection over $S$: set
 \[
 N':= \left\{C=\sum_{\mathrm{finite}}a_iC_i, \hbox{ with }
a_i\in\Z_{\ge0}, \Biggm| \begin{array}{l}
\hbox{where $C_i$ are smooth rational} \\ \hbox{curves on $S$ with
$C^2=-2$} \end{array} \right\}
 \]
and
 \[
 N=\bigl\{C\in N'\setminus \{E_i\} \bigm| (C\cdot E_i)=0
\hbox{ for all }E_i \bigr\}.
 \]
 Then every $C\in N$ induces a reflection $r_C$ on $\coh2{Z}^W$ by 
$r_C(D)=D+( D\cdot C)C$. This is a Hodge isometry.

Set $\Ga=\im\Ga'\to O(S)^+$, where $\Ga'$ is the set of all reflection
$r_C$ for $C\in N$. For this $\Ga$, we have $\Ga \subset O(S)^{++}$ and
the cone $\sA$ is a fundamental domain for the action $\Ga$ on $\sC$.

%%%%%%%%%%%%%%%%%%%%%%%%%%%%%%%%'‹'R,BL(BŽž,BL(BŽå'è-

Recall the weakly polarized global Torelli theorem:

\begin{thm}[Morrison {\cite[Section~4]{M}}] \label{torelli} Let $Y$
and $Y'$ be K3 surfaces with weak polarization, and $\rho\colon Z\to Y$,
$\rho'\colon Z'\to Y'$ their minimal resolutions. Suppose that
$\phi\colon \coh2{Y}^W\to \coh2{Y'}^{W'}$ is an isometry preserving the
Hodge structure such that $\phi(\sA(Y))=\sA(Y')$, which extends to an
isometry $\psi\colon\coh2{Z}\to \coh2{Z'}$. Then there is a unique
isomorphism $\Phi\colon Y\iso Y'$ such that $\Phi^*=\phi$. 
\end{thm}

\begin{lemma}\label{lem4.4}
 \[ 
 \Span{\Ga, P(S)^+}=\Ga \rtimes P(S)^+.
 \]
\end{lemma}
\begin{proof}
 To show that $\Ga$ is a normal subgroup of $\Span{\Ga, P(S)^+}$, we check
that for each $\tau\in P(S)^+$ and $C\in N$, there exists an element
$C'\in N$ such that 
 $\tau^{-1} \circ r_C \circ \tau (z)=r_{C'}$. 
 We have $\tau=f^*$ for some $f\in Q_Z$ then we have:
 \begin{align*}
 \tau^{-1} \circ r_C \circ \tau (z) &=(f^{-1})^* \circ r_C \circ f^*(z) \\
 &=(f^{-1})^*(f^*(z)+f^*(z)\cdot C)C \\
 &=z+(f^*(z)\cdot C)(f^{-1})^*(C) \\
 &=z+(z\cdot (f^{-1})^*C)(f^{-1})^*C.
 \end{align*}
 In addition, $(f^{-1})^*z\in N$ and $f^* \circ g^*=g^* \circ f^*$.
Therefore, we may take $C'=(f^{-1})^*C$.
 
 Assume that $\ga \circ \tau=\ga'\circ\tau'$ for some $\ga, \ga'\in \Ga$
and $\tau, \tau'\in P(S)^+$. Take $\ga''=\ga\circ {\ga'}^{-1}\in \Ga$.
Since $\ga''(\sA)=\sA$, we have $\ga''(\sA)\cap \sA\ne\emptyset$. This
shows that $\ga''=\id$.

 To show that $\Span{\Ga, P(S)^+}\subset \Ga \cdot P(S)^+$, we take
$\tau\in \Span{\Ga, P(S)^+}$, $y\in \sA^{\circ}$ and set $\tau (y)=x$. By
Lemma~\ref{os1.3}, (3), there exists $r_C\in \Ga$ such that $(r_C)^{-1}
\circ \tau(y)=(r_C)^{-1}(x)\in \sA^{\circ}$. Since $S$ is an even
hyperbolic lattice by Lemma~\ref{os1.3}, (1), $(r_C)^{-1} \circ \tau\in
O(S)^+$. Then there exists a Hodge isometry $\rho$ of $\coh2{Z}$ such that
$\rho\rest{S}=(r_C)^{-1} \circ \tau$ and $\rho \circ g^*=g^* \circ \rho$
for all $g\in G$. Recall that $P(S)^+$ is the image of the homomorphism
$Q_Z \ni f\mapsto f^*\rest{S}\in O(S)^{\circ}$. From
Theorem~\ref{torelli}, there exists $f\in \Aut Z$ such that $f=\rho$.
Moreover $f \circ g= g\circ f$ for all $g\in G$ by Theorem~\ref{torelli}
again. Hence $f\in Q_Z$ and $f^*\rest{S}=(r_C)^{-1}\in P(S)^+$. Clearly,
$\Ga\cdot P(S) \subset \Span{\Ga, P(S)^+}$. So we have
$\Span{\Ga,P(S)^+}=\Ga\cdot P(S)^+$. Then we get the result.
\end{proof}

\begin{prop}\label{mainthm}
$\Ga$ is a normal subgroup of $O(S)^+$ and 
 \[
O(S)^+=\Ga \rtimes P(S)^+.
 \]
\end{prop}
\begin{proof}
To prove the first part of the statement, we have only to do with the same
line of Lemma~\ref{lem4.4}.

It is clear that $\Span{\Ga, P(S)^+}\subset O(S)^+$. To show $\Span{\Ga,
P(S)^+}\supset O(S)^+$, take the element $(\tau,\id)\in O(S)\times O(T)$
corresponds to the element of $O(\coh2{Z})$ and this leads to an element of
$O(\coh2{Z}^W)$ by taking restriction. So, applying Theorem~\ref{torelli}
for the pairs $(Y,y)$ and $(Y, (r_C) \circ \tau (y))$, we have the element
of $\Aut Z$ which comes from the element of $\Aut Y$. Then we conclude
that $O(S)^+=\Ga \rtimes P(S)^+$ by Lemma~\ref{lem4.4}.
\end{proof}

\subsubsection*{Proof of Theorem~\ref{th!K3}} Recall that $O(S)^+$ has
finite index in $O(S)^{\circ}$ by Lemma~\ref{os1.3}, (2). Then by
\cite[Chap.~II, p. 116--117]{ash}, there exist a fundamental domain $\De'$
such that $\sC=O(S)^+\cdot \De'$. Changing $\De'$ such as $(\De')^{\circ}
\cap \sA^{\circ}\ne\emptyset$, we may take $\De'$ to be a subset of $\sA$.
(Otherwise, there exists $b\in N$ such that $H_b=\bigr\{x\in S_{\R}
\bigm| (x\cdot b)=0\bigl\}$ which satisfies $(\De)^{\circ}\cap H_b\ne0$.
However, since $r_b(y)=y$ for $y\in (\De)^{\circ}\cap H_b$, we would then
have
$r_b((\De)^{\circ})\ne0$. This is a contradiction.)

By Theorem~\ref{mainthm}, $P(S)^+ \cdot \sA=\sA$, and $\sC=\Ga\cdot\sA$,
we find that $\De'$ is the required fundamental domain.

%%%%%%%%%%%%%%%%%%%%%%%%%%•¶,L#,Hj(B--

\bibliographystyle{amspalpha}

\begin{thebibliography}{99}

\bibitem[AMRT]{ash} A. Ash, D. Mumford, M. Rapport, Y. Tai, \textit{Smooth
compactification of locally symmetric varieties},  Math-Sci.\ Press (1975).

%\bibitem[Bo]{borel} A. Borel, 
%\textit{Introduction aux groupes arithm\'{e}tiques}, 
%Hermann, Paris (1969).

%\bibitem[BPV]{barth} W. Barth, C. Peters, A. Van de Ven,
%\textit{Compact complex surfaces}, 
%Springer-verlag (1984). 

%\bibitem[H]{h} R. Hartshorne, 
%\textit{Algebraic Geometry}, Graduate Texts in Math., 
%Springer-Verlag, Berlin, Heidelberg, and New York (1977).
%\bibitem[Kaw1]{k1} Y. Kawamata, 
%\textit{The cone of curves of algebraic varieties}, Ann. of Math. 
%119 (1984), 603--633. 

\bibitem[Kaw]{k2} Y. Kawamata, \textit{On the cone of divisors of
Calabi--Yau fiber spaces},  Intern.\ J.\ Math.\  8 (1997), 665--687.

%\bibitem[Ko]{ko} S. J. Kov\'{a}cs, 
%\textit{The cone of a K3 surface}, Math. Ann. 300 (1994), 681--691.

\bibitem[M1]{morrison} D. Morrison, \textit{On K3 surfaces with large
Picard number},  Inv.\ Math.\ 75 (1984), 105--121.

\bibitem[M2]{M} D.Morrison, \textit{Some remarks on the moduli of K3
surfaces},  Progress in Math., 29, Birkh\"auser, Boston, Basel,
Stuttgart(1983), 173--259. 

\bibitem[M3]{morrison3} D.Morrison, \textit{Compactifications of Moduli
spaces inspired by Mirror Symmetry},  Journ\'ees de Ge\'om\'eter\'e
Alg\'ebrique d'Orsay, Ast\'erisque 218,  Soc.\ Math.\ France, (1993),
243--291.

\bibitem[O]{type3} K. Oguiso, \textit {On the complete classification of
Calabi--Yau  threefolds of Type} $\mathrm{III}_0$;  in Higher Dimensional
Complex V,  Proc.\ of Int.\ Cong.\ in Trento 1994 (1996), 329--340.

\bibitem[OS]{sakurai} K. Oguiso, J. Sakurai, \textit{Calabi--Yau
threefolds of quotient type}, to appear. preprint math.AG/9909175, 38~pp.

% Oguiso, Keiji An equivariant Torelli theorem for K3 surfaces with
% finite group action and its applications. Free resolutions of coordinate
% rings of projective varieties and related topics (Japanese) (Kyoto,
% 1998). S\=urikaisekikenky\=usho K\=oky\=uroku No. 1078 (1999), 59--63.

\bibitem[S]{shioda} T. Shioda, \textit{The period map of Abelian
surfaces}, J.\ Fac.\ Sci.\ Univ.\ Tokyo Sect. \  IA25 (1978), 47--59.

\bibitem[St]{sterk} H. Sterk, \textit{Finiteness results for algebraic K3
surfaces}, Math.\ Z.\ 189 (1985), 507--513.

\bibitem[Zh]{Zhang} D.-Q. Zhang, \textit{Logarithmic Enriques surfaces},
J.\ Math.\ Kyoto Univ.\ 31 (1991), 419--466.

\end{thebibliography}

\end{document}